\documentclass[ampa]{svjour}
\usepackage{graphics,amsmath,amssymb,amsfonts,amsopn}
\renewcommand{\ge}{\varepsilon}
\newcommand{\de}{\partial}
\newcommand{\R}{\mathbb{R}}

\begin{document}
\author{Simone Secchi \and Enrico Serra}
\title{Symmetry breaking results for problems with exponential growth
in the unit disk\thanks{This research was supported by MURST
``Variational Methods
and Nonlinear Differential Equations''}}
\titlerunning{Symmetry breaking results}
\institute{S.~Secchi\at Dipartimento di Matematica, Universit\`a di
  Milano\\ \email{secchi@mat.unimi.it} \and E.~Serra\at Dipartimento di
  Matematica, Universit\`a di Milano \\
  \email{enrico.serra@mat.unimi.it}}
\date{\today}
 \maketitle
%
\begin{abstract}
We investigate some asymptotic properties of extrema $u_\alpha$ to the
two-dimensional variational problem
\begin{equation*}
\sup_{\substack{u\in H_0^1(B) \\ \| u\|=1}} \int_B
\left( e^{\gamma u^2} -1 \right) |x|^\alpha \, dx
\end{equation*}
as $\alpha \to +\infty$. Here $B$ is the unit disk of $\R^2$ and
$0<\gamma\leq 4\pi$ is a given parameter. We prove that in a certain
range of $\gamma$'s, the maximizers are not radial for $\alpha$ large.
\end{abstract}

\section{Introduction}
\label{intro}

Let $0<\gamma\leq 4\pi$ be a given number. We consider the
maximization problem
\begin{equation} \label{P}
S(\alpha,\gamma)=\sup_{\substack{u\in H_0^1(B) \\ \| u\|=1}}
\int_B
\left( e^{\gamma u^2} -1 \right) |x|^\alpha \, dx,
\end{equation}
where $B=\left\{ x \in \R^2 : |x|<1 \right\}$, $H_0^1(B)$ is the usual
Sobolev space endowed with the Dirichlet norm $\|u\|=\left( \int_B
|\nabla u|^2 \, dx \right)^{1/2}$, and $\alpha >0$. It is readily seen
that any maximizer of \eqref{P} must satisfy (weakly) the elliptic differential
equation
\begin{equation} \label{eulero}
-\Delta u = \lambda |x|^\alpha u e^{\gamma u^2}
\end{equation}
with a Lagrange multiplier $\lambda$ given by
\begin{equation}
\label{lambda}
\lambda = \frac{1}{\int_B u^2e^{\gamma u^2}|x|^\alpha \, dx}
\end{equation}
Standard regularity
theory shows that any weak solution of (\ref{eulero}) is classical.
Moreover, if $u$ is a maximizer of problem \eqref{P}, then so is $|u|$. Hence we can
work with nonnegative functions. We will use
freely these facts.

Our problem can be seen as a natural two--dimensional extension of the
H\'{e}non--type problem
\begin{equation}
\label{henon1}
\sup_{\substack{u\in H_0^1(B) \\ \| u\|=1}} \int_B |u|^p |x|^\alpha \, dx
\end{equation}
in $\R^n$ with $n \geq 3$ and $1<p<2^*$. Indeed, by the
Trudinger--Moser inequality (see \cite{moser,ruf,trudinger}),
\begin{equation} \label{trud}
\sup_{\substack{u\in H_0^1(B) \\ \| u\|=1}} \int_B \left( e^{\gamma
u^2} -1 \right) \, dx
\begin{cases}
<\infty &\text{if $\gamma\leq 4\pi$} \\ =+\infty &\text{otherwise,}
\end{cases}
\end{equation}
the growth $\exp(4 \pi |\cdot|^2)$ in 2D corresponds
(with relevant differences, though) to the
critical growth $| \cdot|^{2^*}$ in dimension $n \geq 3$. This can be
made precise by introducing the class of Orlicz spaces, but we shall
not go into the details. We refer the interested reader to \cite{a,ct}.

Recently, Smets \textit{et al.} (\cite{ssw}) studied the symmetry of
minimizers to the problem
\begin{equation} \label{s1}
S_{\alpha,p} = \inf_{\substack{u\in H_0^1(B) \\ u \neq 0}}
\frac{\int_B |\nabla u|^2 \, dx}{\left( \int_B |x|^\alpha |u|^p \,
dx\right)^{2/p}}\;,
\end{equation}
namely problem (\ref{henon1}), with
\begin{eqnarray*}
&&2<p<+\infty \quad \hbox{in dimension 2} \\
&&2<p<2^* \quad \hbox{in higher dimension.}
\end{eqnarray*}
Since the quotient in (\ref{s1}) is invariant under rotations, it is natural
to set up the same minimization problem in the space of \textit{radial}
functions
$H_{0,{\rm rad}}^1(B)$:
\begin{equation} \label{s2}
S_{\alpha,p}^{\rm rad} = \inf_{\substack{u\in H_{0,{\rm rad}}^1(B) \\
    u \neq 0}} \frac{\int_B |\nabla u|^2 \, dx}{\left( \int_B
    |x|^\alpha |u|^p \, dx\right)^{2/p}}.
\end{equation}
The set $B$ being bounded, both problems \eqref{s1} and \eqref{s2}
are compact, and are thus solved by functions $u_\alpha$ and $v_\alpha$. A
very interesting \textit{symmetry--breaking} result contained in
\cite{ssw} is the following.
\begin{theorem}[Smets, Su, Willem]
  Assume the dimension of the space is greater than or equal to 2. For
  every $p\in (2,2^*)$ ($p>2$ in 2D), there exists $\alpha^*>0$ such
  that no minimizer of \eqref{s1} is radial provided that $\alpha >
  \alpha^*$. In particular,
\[
S_{\alpha,p} < S_{\alpha,p}^{\rm rad} \quad \hbox{for all $\alpha$
sufficiently large}.
\]
\end{theorem}
This result has generated a line of research on H\'enon--type equations.
For example, it shows in particular that for a certain
set of parameter
values, the (H\'enon) equation associated to (\ref{s1}) admits the
coexistence of  radial and  nonradial positive solutions. Since a radial
solution always exists, by a result of Ni, \cite{Ni}, for $\alpha <
2^* + \frac{2\alpha}{n-2}$, a similar phenomenon can be expected
also for critical and supercritical growths. Results in this direction
have been obtained in \cite{se} and \cite{bs}. See also \cite{bw} and
\cite{SW} for asymptotic analysis of ground states and other symmetry results.
\medskip

The symmetry breaking problem for exponential nonlinearities in the unit disk,
on the contrary, seems to have been much less studied.
Very recently, Calanchi and Terraneo (see \cite{ct}) proved some
results about the existence of non--radial maximizers for the
variational problem
\begin{equation*}
T_{\alpha,p,\gamma} = \sup_{\substack{u\in H_0^1(B) \\ \|u\| \leq 1}}
\int_B |x|^\alpha \left( e^{p|u|^\gamma} -1 -p|u|^\gamma \right)\, dx
\end{equation*}
where $\alpha>0$, $p>0$ and $1<\gamma \leq 2$ when $\alpha \to
+\infty$.  We observe that the functional to be maximized in
$T_{\alpha,p,\gamma}$ contains an extra term with respect to our
$S(\alpha,\gamma)$.

\bigskip

In this paper we present some results about
symmetry of solutions to \eqref{P}, though we are not able to
cover the whole range $(0,4\pi]$ of the parameter $\gamma$. The
main difficulty is that, unlike \eqref{s1}, our problem \eqref{P}
is not homogeneous with respect to $u$. As a consequence, we
cannot replace \eqref{P} with a more familiar ``Rayleigh''
quotient.

We consider problem \eqref{P} and
its radial companion
\begin{equation} \label{comp}
S^{\rm rad}(\alpha,\gamma) = \sup_{\substack{u\in H_{0,{\rm rad}}^1(B)
    \\
\|u\|=1}} \int_B \left( e^{\gamma u^2} -1 \right) |x|^\alpha \, dx;
\end{equation}
since $H_{0,{\rm rad}}^1(B) \subset H_0^1(B)$, it is
clear that
\begin{equation}
\label{ineq}
S(\alpha,\gamma) \geq S^{\rm rad}(\alpha,\gamma).
\end{equation}
Our main concern is to investigate if and when the {\em strict} inequality
takes place.
By standard arguments (see Section 2), both
$S(\alpha,\gamma)$ and $S^{\rm rad}(\alpha,\gamma)$ are attained
for $\gamma \in (0,4\pi)$ (this interval is considerably larger for
the radial case, see \cite{ct}).

We first obtain an asymptotic profile type result for the maximizers of
(\ref{comp}), as $\alpha \to \infty$. This result is essential in order
to carry out the proof of the main symmetry breaking theorem, and we believe
that it is interesting in its own.

In the statements that follow we denote by $\lambda_1$ the first eigenvalue
of $-\Delta$ on $H_0^1(B)$, and by $\varphi_1$ the corresponding (positive)
eigenfunction, normalized by $||\varphi_1||=1$.

\begin{theorem}
\label{profile}
Let $\gamma \in (0,4\pi]$. For every $\alpha >0$,
let $u_\alpha = u_\alpha(|x|)$ be a maximizer for $S^{\rm rad}
(\alpha,\gamma)$. Then
$$
\lim_{\alpha \to +\infty}{\scriptstyle \sqrt{\frac{\alpha+2}{2}}}
\, u_\alpha (|x|^{\frac{2}{\alpha+2}}) = \varphi_1 \quad
\hbox{strongly in $H_0^1(B)$}.
$$
\end{theorem}

We now state the main result of the paper.

\begin{theorem}
\label{main} There exists $\gamma^* \in [0,4\pi)$ such that for
every $\gamma \in (\gamma^*,4\pi)$ no maximizer for $S(\alpha,\gamma)$
is radial provided $\alpha$ is large enough. Moreover
\begin{equation}
\label{gammastar}
\gamma^* \leq \frac{\pi \varphi_1(0)^2}{\lambda_1 \int_B \varphi_1^4 \, dx}.
\end{equation}
\end{theorem}
Of course the upper bound for $\gamma^*$ appearing in the right--hand--side of
(\ref{gammastar}) is strictly smaller than $4\pi$.
We do not know if $\gamma^* =0$; this is one of the interesting open problems
connected to $S(\alpha,\gamma)$ and should be the object of further research.
\medskip

The paper is structured as follows. In Section 2 we obtain the asymptotic
description of radial maximizers, while a similar result, perturbative
in nature, is given in Section 3 for non radial maximizers. Section 4
is devoted to the proof of the main result, Theorem \ref{main}.

\section{Asymptotic behavior of radial maximizers}

In this section we give a precise description of the asymptotic
behavior of maximizers of problem~\eqref{comp} as $\alpha \to
+\infty$.

To begin with, we fix some notation that we will use throughout the paper.
We introduce the \textit{variational functional}
\begin{equation} \label{funct}
I(u)= \int_B \left( e^{\gamma u^2}-1 \right) |x|^\alpha \, dx
\end{equation}
which acts formally in the same way both on $H_0^1(B)$ and on
$H_{0,{\rm rad}}^1(B)$. For the sake of simplicity we suppress the
dependence of $I$ on $\alpha$ and $\gamma$.
The first (positive) eigenfunction of the
Laplace operator $-\Delta$ on $H_0^1(B)$ will be denoted by
$\varphi_1$, normalized by $||\varphi_1|| = 1$, and the corresponding
eigenvalue by $\lambda_1$.

Throughout the paper, we will make use of polar coordinates in $\R^2$,
namely $x=(\rho \cos \theta, \rho\sin \theta)$ with $\rho \geq 0$ and $\theta
\in [0,2\pi)$. With a slight abuse of notation, we will write
$u(\rho,\theta)=u(x)=u(\rho \cos \theta, \rho\sin \theta)$ for a given function
$u$ on $\R^2$, and, likewise, $u(x) = u(|x|) = u(\rho)$ for a radial
function. For further use, we state the variational
problem \eqref{P} in polar coordinates.
Set
\begin{equation} \label{eps}
\ge = \frac{2}{\alpha +2} \raise 2pt \hbox{.}
\end{equation}
For any smooth (or $H_0^1$) function $u$ on $B$, define the new function
\begin{equation} \label{v}
v(\rho,\theta):= \frac{1}{\sqrt \ge} u(\rho^\ge,\theta)
\end{equation}
expressed in polar coordinates. Observe that
\begin{align*}
v_\rho &= \sqrt{\ge} \rho^{\ge -1} u_\rho (\rho^\ge,\theta) \\
v_\theta &=
\frac{1}{\sqrt{\ge}} u_\theta (\rho^\ge,\theta),
\end{align*}
so that, if $t=\rho^\ge$,
\begin{equation} \label{vinc:pol}
\int_0^1 \int_0^{2\pi} \left(v_t^2 +
\frac{\ge^2}{t^2} v_\theta^2 \right) t\, dt\, d\theta=1
\end{equation}
whenever $\int_B |\nabla u|^2 \, dx=1$. In the variables $(t,\theta)$,
the variational functional \eqref{funct} reads
\begin{equation*}
I(v) = \ge \int_0^1 \!\int_0^{2\pi} \left( e^{\gamma \ge v^2} -1
\right) t \, dt\, d\theta.
\end{equation*}
Finally, the original problem \eqref{P} can be written by means
of \eqref{eps}, \eqref{v} as
\begin{multline} \label{P:polar}
S(\alpha,\gamma)= \sup \Bigg\{ \ge \int_0^{1} \int_0^{2\pi} \left(
e^{\ge\gamma v^2} -1 \right) t\, dt \, d\theta : \\ \int_0^1 \int_0^{2\pi} \left(v_t^2 +
\frac{\ge^2}{t^2} v_\theta^2 \right) t\, dt\, d\theta=1 \Bigg\}.
\end{multline}

\begin{remark}
  We stress that in the new variables the weight $|x|^\alpha$
  disappears from the functional and the parameter $\varepsilon =
  \frac{2}{\alpha+2}$ appears both in the exponent and
  in front of $|\de v / \de \theta|^2$. Notice that if $u$ is radial,
  then
\begin{equation}
 \label{eq:rad1}
\int_B |\nabla u|^2 \, dx = 2\pi \int_0^1 v_t^2 t \, dt
= \int_B |\nabla v|^2 \, dx
\end{equation}
and
\begin{equation}
\label{eq:rad2}
\int_B \left( e^{\gamma u^2} -1 \right)|x|^\alpha \, dx = 2\pi\ge \int_0^1
\left( e^{\ge \gamma v^2}-1 \right)t\,dt
=\varepsilon\int_B \left( e^{\varepsilon\gamma v^2} -1 \right) \, dx,
\end{equation}
so that
\begin{equation} \label{P:radial}
S^{\rm rad}(\alpha,\gamma)= \sup \Bigg\{2\pi \ge \int_0^{1}  \left(
e^{\ge\gamma v^2} -1 \right) t\, dt  :  2\pi\int_0^1 v_t^2  t\, dt=1 \Bigg\}.
\end{equation}
\end{remark}
First of all, we deal with the existence of maximizers to
\eqref{P} and \eqref{comp}.

\begin{proposition} \label{prop1}
There exist a solution $u_\alpha \in
H_0^1(B)$ and a solution $u_\alpha^r \in H_{0,{\rm rad}}^1 (B)$ (also
called $u_\ge$ and $u_\ge^r$ via the change of parameter \eqref{eps})
to problems~\eqref{P} and \eqref{comp} respectively, provided
$0<\gamma < 4\pi$.

\noindent When $\gamma=4\pi$, problem \eqref{comp} has a solution in $H_{0,{\rm
rad}}^1 (B)$.
\end{proposition}
\begin{proof}
We only give some details, since the argument can be recovered from
the existing literature.

First of all we notice that
\[
S(\alpha,\gamma) = \sup_{\substack{u\in H_0^1(B) \\ \| u\| \leq 1}}
\int_B
\left( e^{\gamma u^2} -1 \right) |x|^\alpha \, dx,
\]
and the same for $S^{\rm rad}(\alpha,\gamma)$.
In the \textit{subcritical case} $\gamma <
4\pi$, the proof is almost trivial.  Indeed, let $\{u_n\}$ be a
maximizing sequence for $S(\alpha,\gamma)$ (or for
$S^{\rm rad}(\alpha,\gamma)$), with $\|u_n\| \leq 1$. By
the Sobolev embedding theorem, we can assume without loss of
generality that (up to a subsequence) $u_n \rightharpoonup u$,
weakly in $H_0^1$ and $u_n
\to u$ a.e. and strongly in $L^q(B)$ for any finite $q \geq 1$. In
particular, $\|u\| \leq 1$. Then, thanks to Lemma 2.1 of
\cite{dfmr}, we have
\[
S(\alpha,\gamma)=
\lim_{n \to \infty} \int_B \left( e^{\gamma u_n^2} -1 \right)
|x|^\alpha\, dx = \int_B \left( e^{\gamma u^2} -1 \right)
|x|^\alpha\, dx.
\]
This shows that $u \neq 0$, and that $u$ is a maximizer of \eqref{P}.

The \textit{critical case} $\gamma=4\pi$ for $S^{\rm rad}(\alpha,\gamma)$ is
slightly different. Indeed, equation \eqref{eq:rad2} shows that
problem \eqref{comp} is still ``subcritical'', provided that $\ge
\gamma < 4\pi$, i.e. $\gamma < 4\pi + 2\pi \alpha$. Therefore,
standard arguments prove that $S^{\rm rad}(\alpha,4\pi)$ is
actually attained by a radial function. See also the remark at the
end of section 3 in~\cite{ct}. \qed\end{proof}

\begin{remark}
It does not seem to be known whether
$S(\alpha,4\pi)$ is attained. For the ``unweighted case''
$\alpha=0$ this is a  celebrated result due to Carleson and Chang
\cite{cc}. Unfortunately, it does not seem possible to modify
their proof so as to take into account the weight $|x|^\alpha$.
This is an interesting open problem.
\end{remark}

We now begin the study of the asymptotic behavior of the radial maximizers.

Take a radial function $v$, compactly supported in $B$, with
$|| v|| =1$. Formally,
\begin{align}
\label{expans}
e^{\ge\gamma v^2} -1 &= \ge \gamma v^2 + \frac{1}{2!} \ge^2 \gamma^2
v^4 + \frac{1}{3!} \ge^3 \gamma^3 v^6 + \dots \nonumber \\ &=: \ge \gamma v^2 +
R_\ge (v)
\end{align}

\begin{lemma}
As $\ge \to 0$,
\begin{equation} \label{lem1.1}
\ge \int_0^1 \left( e^{\ge \gamma v^2}-1 \right) t \, dt\, d\theta
= \ge^2 \gamma \int_0^1 v^2 t \, dt + O(\ge^3),
\end{equation}
uniformly for $||v||=1$.
\end{lemma}

\begin{proof}
Equation \eqref{lem1.1} is equivalent, via (\ref{expans}), to
\begin{equation}
\int_0^1 R_\ge (v)t \, dt = O(\ge^2).
\end{equation}
Now, $R_\ge (v) = \sum_{k=2}^\infty \frac{1}{k!} (\ge \gamma)^k
v^{2k}$. Fix any index $k \geq 2$: for every $t\in [0,1]$, we have
\begin{align*}
v(t) &= v(t)-v(1)=\int_1^t v'(s) \, ds
\le \left(\int_t^1 |v'(s)|^2 s\,ds\right)^{1/2}
\left(\int_t^1 \frac1s\,ds\right)^{1/2}
\\ &\leq \frac{1}{\sqrt{2\pi}}\left( \int_B
|\nabla v|^2 \, dx \right)^{1/2} \left( \int_1^t \frac{ds}{s}
\right)^{1/2} = \frac{1}{\sqrt{2\pi}} \left( -\log t \right)^{1/2}.
\end{align*}
As a consequence, $|v(t)|^{2k} \leq \frac{1}{(2\pi)^k} \left( -\log t
\right)^k$. We multiply by $t$, integrate this inequality
over $[0,1]$ and find
\[
\int_0^1 |v(t)|^{2k} t \, dt \leq \frac{1}{(2\pi)^k} \int_0^1 \left( -\log t \right)^k t \, dt.
\]
The change of variable $t=\exp(-x/2)$ yields immediately
\begin{align}
\int_0^1 |v(t)|^{2k} t \, dt &\leq \frac{1}{(2\pi)^k} \int_0^\infty
\frac12 e^{-x/2} e^{-x/2} \left( \frac{x}{2} \right)^k \, dx \notag \\
& = \frac{1}{2^{k+1}} \frac{1}{(2\pi)^k} \int_0^\infty x^k e^{-x} \,
dx \notag \\ &= \frac{1}{2^{k+1}} \frac{1}{(2\pi)^k} \Gamma (k+1) =
\frac{1}{2^{k+1}} \frac{k!}{(2\pi)^k}.
\end{align}
The Monotone Convergence Theorem implies that we can switch the
summation over $k$ with the integration over $[0,1]$, so that
\begin{align*}
\int_0^1 R_\ge (v)t \, dt &= \sum_{k=2}^\infty \int_0^1 \frac{1}{k!}
(\ge \gamma)^k |v(t)|^{2k} t \, dt \leq \sum_{k=2}^\infty
\frac{(\ge\gamma)^k}{k!} \frac{1}{(2\pi)^k} \frac{k!}{2^{k+1}} \\
&= \sum_{k=2}^\infty \frac{(\ge\gamma)^k}{2(4\pi)^k} =
\frac12 \sum_{k=2}^\infty \left( \frac{\ge\gamma}{4\pi} \right)^k
=\frac{\gamma \ge^2}{8\pi (4\pi - \gamma \ge)} = O(\ge^2)
\end{align*}
This completes the proof.
\qed\end{proof}
We now establish the asymptotic behavior of a sequence $u_\alpha$ of maximizers
of $S^{\rm rad}(\alpha,\gamma)$ as $\alpha\to +\infty$.
For notational convenience, in the statement of the result we denote
this sequence by $u_\varepsilon$, keeping in mind that $\alpha$ and
$\varepsilon$ are linked by (\ref{eps}).

\begin{theorem} \label{th2} Let $\gamma\in (0,4\pi]$ and
let $u_\varepsilon \in H_{0,{\rm rad}}^1(B)$ be a maximizer of
$S^{\rm rad}(\alpha,\gamma)$.
Then
\begin{equation} \label{converge}
\lim_{\ge \to 0}\frac{1}{\sqrt\varepsilon} u_\varepsilon(|x|^\varepsilon) =
\varphi_1 \quad \hbox{strongly in $H_0^1(B)$}.
\end{equation}
\end{theorem}

\begin{proof}

We set $v_\varepsilon(t) = \frac{1}{\sqrt\ge}u_\ge (t^\ge)$. Clearly
$v_\varepsilon$ is a maximizer of problem
\eqref{P:radial}. In particular, $\|v_\ge \|=1$ for all $\varepsilon$,
so that the set $\{v_\ge\}_{\ge}$ is bounded in
$H_0^1(B)$; therefore some subsequence, which we still term
$\{v_\ge\}_{\ge}$, converges weakly to some $v$ in $H_0^1(B)$,
and strongly in $L^q(B)$ for all finite $q\geq 1$ as $\ge \to 0$.
Since $v_\varepsilon$ is a maximizer,
for every radial function $\psi\in H_0^1(B)$ satisfying $||\psi||=1$ we have
\begin{align*}
2\pi\ge  \int_0^1 \left( e^{\ge \gamma v_\ge^2} -1 \right)t
\, dt &\geq 2\pi\ge  \int_0^1 \left( e^{\ge \gamma
\psi^2} -1 \right)t \, dt \\ &\geq 2 \pi \gamma \ge^2
\int_0^1 \psi^2 t \, dt.
\end{align*}
Hence by Lemma 1
\[
2 \pi \gamma \ge^2 \int_0^1 \psi^2 t \, dt \leq 2\pi\gamma \ge^2
\int_0^1 v_\ge^2 t \, dt + O(\ge^3).
\]
Dividing by $\gamma\varepsilon^2$ and letting $\ge \to 0$ we obtain
\begin{equation} \label{psi}
2\pi\int_0^1 \psi^2 t \, dt \leq 2\pi\int_0^1 v^2 t \, dt,
\end{equation}
namely
\begin{equation}
\label{ineql}
\int_B \psi^2\,dx \le \int_B v^2\,dx.
\end{equation}
If we now maximize over those $\psi\in H_0^1(B)$ satisfying
$\| \psi\|=1$ we see that
\begin{equation} \label{v-minore}
\frac{1}{\lambda_1} \leq  \int_B v^2 \, dx.
\end{equation}
This shows that $v \neq 0$. By a standard semicontinuity
argument, we have $|| v||  \leq 1$. Therefore
$$
\lambda_1 \leq \frac{\int_B |\nabla v|^2 \, dx}{\int_B v^2 \, dx} \leq
\frac{1}{\int_B v^2 \, dx} \leq \lambda_1,
$$
which shows that
$$
\int_B |\nabla v|^2 \, dx= 1 \qquad\hbox{and}\qquad \int_B v^2 \, dx =
\frac{1}{\lambda_1}.
$$
Since $\lambda_1$ is a simple eigenvalue, this means that
$v=\varphi_1$ and
\begin{equation}
\int_B v_\ge^2 \, dx  \to \int_B \varphi_1^2 \, dx,
\quad \int_B |\nabla v_\ge|^2 \, dx \to \int_B |\nabla \varphi_1|^{2}.
\end{equation}
This, together with the weak convergence $v_\ge \rightharpoonup
v=\varphi_1$, shows that $v_\ge \to \varphi_1$ strongly in
$H_0^1(B)$.
\qed\end{proof}
We complete the description of the asymptotics
of problem (\ref{comp}) with the behavior of the levels.

\begin{proposition} \label{prop-rad}
Let $\gamma \in(0,4\pi]$. As $\ge \to 0$, i.e. as $\alpha \to
+\infty$, we have
\begin{equation} \label{eq:rad-as}
S^{\rm rad} (\alpha,\gamma) = \frac{\gamma}{\lambda_1} \ge^2 + o(\ge^2).
\end{equation}
\end{proposition}

\begin{proof}

Inserting $v=v_\ge$, as defined in the proof of the previous theorem,
in \eqref{lem1.1}, we obtain
\begin{multline*}
2\pi\ge \int_0^1 \left( e^{\gamma \ge v_\ge^2}-1 \right) t
\, dt  = 2\pi \gamma \ge^2 \int_0^1 v_\ge^2 t \, dt
+O(\ge^3) \\ =\gamma \ge^2 \left( \frac{1}{\lambda_1} +o(1) \right)
+O(\ge^3) = \frac{\gamma}{\lambda_1} \ge^2 + o(\ge^2).
\end{multline*}
\qed
\end{proof}

\section{Asymptotic estimates for non--radial maximizers}

In the previous section we have proved that $S^{\rm rad}
(\alpha,\gamma) \approx \frac{\gamma}{\lambda_1} \ge^2$ when $\ge =
2/(\alpha+2) \to 0$. We now provide a similar estimate for
$S(\alpha,\gamma)$, and show that solutions to \eqref{P} are never
radial, provided $\alpha$ is large and $\gamma \approx 4\pi$.

We begin with a lemma which estimates $S(\alpha,\gamma)$ in terms
of the functional without weight.

\begin{lemma}
\label{transf}
Let $\gamma \in (0,4\pi]$ and $\alpha>0$. Setting as usual $\varepsilon =
2/(\alpha+2)$, we have
\begin{equation}
\label{eq:transf}
S(\alpha,\gamma) > \frac{\varepsilon^2}{4}\sup_{||v|| = 1}
\int_B \left(e^{\gamma v^2} - 1\right)\, dx.
\end{equation}
\end{lemma}

\begin{proof} We denote by $p$ be the point $(-\frac12,0)\in B$ and we take
a function $\psi \in H_0^1(B_{1/2}(p))$. Notice that in polar
coordinates, the function $\psi$ vanishes on
$\partial([0,1]\times[0,2\pi))$. We extend then $\psi$ by zero
outside $[0,1]\times[0,2\pi)$ and we still call $\psi$ this
extension. Then it makes sense to define, for every $\varepsilon
\in(0,1)$, a function $u = u_\varepsilon : [0,1]\times[0,2\pi)\to
\R$ by
$$
u(t,\phi) = \psi(t^{1/\varepsilon},\phi/\varepsilon).
$$
Obviously the function $u$ is not radial.
Setting $\rho = t^{1/\varepsilon}$ and $\theta = \phi/\varepsilon$ we see that
\begin{align*}
&\int_B|\nabla u|^2\,dx = \int_0^1  \int_0^{2\pi}
(u_t^2 +\frac{1}{t^2}u_\phi^2)t\,dt\,d\phi \\  =
&\int_0^1\int_0^{2\pi} \left(\frac{t^{2/\varepsilon - 2}}{\varepsilon^2}
\psi_\rho^2(t^{1/\varepsilon},\phi/\varepsilon) + \frac{1}{\varepsilon^2 t^2}
\psi_\theta^2(t^{1/\varepsilon},\phi/\varepsilon)\right)t\,dt\,d\phi\\
 = &\int_0^1\int_0^{2\pi} \left( \psi_\rho^2 + \frac{1}{\rho^2}\psi_\theta^2\right)
\rho\,d\rho d\theta = \int_B |\nabla \psi|^2\, dx =
\int_{B_{1/2}(p)} |\nabla \psi|^2\, dx
\end{align*}
and

\begin{align*}
&\int_B\left(e^{\gamma u^2} -1\right)|x|^\alpha\,dx =
\int_0^1\int_0^{2\pi} \left(e^{\gamma u^2} -1\right)t^{\alpha+1}\,dt \, d\phi \\
=& \int_0^1\int_0^{2\pi} \left(e^{\gamma
\psi^2(t^{1/\varepsilon},\phi/\varepsilon)}
-1\right)t^{\alpha+1}\,dt \, d\phi =
\varepsilon^2\int_0^1\int_0^{2\pi} \left(e^{\gamma \psi^2}
-1\right)
\rho\,d\rho \, d\theta \\
= &\varepsilon^2\int_B  \left(e^{\gamma \psi^2} -1\right)\,dx =
\varepsilon^2\int_{B_{1/2}(p)} \left(e^{\gamma \psi^2} -1\right)\,dx.
\end{align*}
Therefore we can say that

\begin{align}
\label{medium}
S(\alpha,\gamma) >
\sup\Bigg\{\varepsilon^2\int_{B_{1/2}(p)} &\left(e^{\gamma \psi^2} -1\right)\,dx
\;: \nonumber\\
&\psi \in H_0^1(B_{1/2}(p)),\; \int_{B_{1/2}(p)}|\nabla \psi|^2\,dx = 1 \Bigg\}.
\end{align}
Next we define $v\in H_0^1(B)$ as $v(x) = \psi(x/2 +p)$. Obviously
$$
\int_B |\nabla v|^2 \,dx = \int_{B_{1/2}(p)} |\nabla \psi|^2 \,dx,
$$
while
$$
\int_B \left(e^{\gamma v^2} -1\right)\,dx = 4
\int_{B_{1/2}(p)} \left(e^{\gamma \psi^2} -1\right)\,dx.
$$
This means, by (\ref{medium}), that
\begin{equation}
\label{final}
S(\alpha,\gamma) > \frac{\varepsilon^2}{4}\sup_{||v|| = 1}
\int_B \left(e^{\gamma v^2} - 1\right)\, dx,
\end{equation}
and the proof is complete. \qed
\end{proof}
We are now ready to state the main result of this section.
\begin{theorem} \label{prop4}
There exists $\gamma^* \in (0,4\pi]$ such that, for all $\gamma \in
(\gamma^*, 4\pi)$,
\begin{equation} \label{eq:main}
S(\alpha,\gamma) > S^{\rm rad}(\alpha,\gamma),
\end{equation}
provided $\alpha$ is large enough.
\end{theorem}

\begin{proof}
By the results of Section 2 we know that
$$
S^{\rm rad} (\alpha,\gamma) = \frac{\gamma}{\lambda_1} \ge^2 + o(\varepsilon^2)
$$
as $\varepsilon \to 0$. In view of (\ref{final}) the proof is done if we show that
\begin{equation}
\label{toprove}
\frac{1}{4}\sup_{||v|| = 1}
\int_B \left(e^{\gamma v^2} - 1\right)\, dx > \frac{\gamma}{\lambda_1}.
\end{equation}
The value in the left--hand side of \eqref{toprove}, which is
attained by the results in \cite{cc}, is unknown. We are going to
estimate it using the same function that appears in \cite{cc}.

From now on, we assume that $v$ is a radial function. This is
natural, since the supremum in (\ref{toprove}) is attained by a
radial function. If $v$ is radial, it is convenient to introduce
the function $w: [0,+\infty) \to \R$ defined by
\begin{equation*}
w(t):=\sqrt{4\pi} v(e^{-t}).
\end{equation*}
A straightforward computation shows that
\begin{equation}
\int_B |\nabla v|^2\, dx = \int_0^\infty |w'(t)|^2 \, dt
\end{equation}
and
\begin{equation*}
\int_B \left( e^{\gamma v^2}-1 \right)\, dx =
\pi \int_0^\infty \left( e^{\frac{\gamma}{4\pi}w^2(t)}-1 \right) e^{-t} \, dt.
\end{equation*}
Since the statement of the Theorem is perturbative in nature with
respect to $\gamma$, and everything depends continuously on
$\gamma$, to complete the proof we can assume $\gamma = 4\pi$.
Explicitly, we focus on the problem
\begin{equation*}
\max_{\int_0^\infty
|w'|^2 \, dt=1}  \pi \int_0^\infty e^{w^2-t}\, dt -\pi
\end{equation*}
Take then $w \colon [0,+\infty) \to \R$ to be
\begin{equation*}
w(t)=
\begin{cases}
\frac{1}{2} t &\text{if $0 \leq t \leq 2$} \\ \sqrt{t-1} &\text{if $2
\leq t \leq 1+e^2$} \\ e &\text{if $t \geq 1+e^2$}.
\end{cases}
\end{equation*}
This is the function that already appears in \cite{cc}. By direct inspection,
$$
\int_0^\infty
|w'|^2 \, dt=1\qquad\hbox{and}\qquad
\int_0^\infty e^{w^2-t}\, dt = \frac{2}{e} \int_0^1 e^{t^2}\, dt + e.
$$
Hence
\begin{equation}
\max_{\int_0^\infty |w'|^2 \, dt=1} \pi \int_0^\infty e^{w^2-t}\, dt -\pi
> \pi \left( \frac{2}{e} \int_0^1
e^{t^2}\, dt + e \right) - \pi.
\end{equation}
If we show that
\begin{equation} \label{finish-prop}
\frac{1}{4} \left( \frac{2\pi}{e} \int_0^1 e^{t^2}\, dt + e\pi
-\pi \right) > \frac{4\pi}{\lambda_1},
\end{equation}
then also (\ref{toprove}) will be satisfied, by continuity, for $\gamma$
close enough to $4\pi$,
and the proof will be finished.
We thus check that
$$
\frac{2}{e} \int_0^1 e^{t^2}\, dt + e-1 > \frac{16}{\lambda_1}.
$$
From the characterization of $\lambda_1$ as a zero of the Bessel
function $J_0$ (see \cite{ch,web}), we have the approximated value
$\lambda_1 \approx 5.783$.  If we estimate $\int_0^1 e^{t^2}\, dt
$ by expanding the integrand in power series, and taking into
account only the first three terms, we get easily that $\int_0^1
e^{t^2}\, dt > 1.453$.\footnote{A good approximation provided by
Maple\circledR\ for the integral is $\int_0^1 \exp (t^2)\, dt
\approx 1.462651746$.} Therefore
$$
\frac{2}{e} \int_0^1 e^{t^2}\, dt + e-1 >
\frac{2.906}{e}+e-1 \approx 2.787 >  2.767  \approx \frac{16}{\lambda_1}.
$$
\qed\end{proof}
Although we do not know whether problem \eqref{P} admits a solution
in the critical case $\gamma=4\pi$,
the previous proof gives the following \textit{a
priori} information.

\begin{corollary} \label{cor:1}
If $S(\alpha,4\pi)$ is attained by some function $u$, then $u$ cannot be radial.
\end{corollary}

\section{A nonperturbative estimate for symmetry--breaking}

We have seen in Theorem~\ref{prop4} that solutions to \eqref{P}
are non-radial whenever $\gamma$ is
close to $4\pi$ and $\alpha$ is large (depending on $\gamma$).
In this final section we present a similar
result, whose nature is no longer perturbative with respect to
$\gamma$. The technique of the proof is rather different, and
resembles that of Theorem 2.1 in \cite{ssw}.

For clarity purposes,
we introduce an auxiliary map $N$, defined by
\begin{equation} \label{N}
N(u)=\frac{u^2}{|| u||^2 } = \frac{u^2}{\int_B |\nabla u|^2\,dx }
\quad \hbox{for all $u\in
H_0^1(B)\setminus\{0\}$}.
\end{equation}
and a measure $\mu_\alpha$ on Borel subsets $E$ of $\R^2$ by
$$
\mu_\alpha (E) = \int_E |x|^\alpha \, dx.
$$
It follows from straightforward arguments that problem \eqref{P} is
equivalent to the maximization of the ``free'' functional
\[
F(u)=\int_B \left(e^{\gamma N(u)}-1\right)\, d\mu_\alpha
\]
on $H_0^1(B)\setminus \{0\}$. The use of this functional allows us
to embed some homogeneity in the problem, which will be very useful
in the computations below.

In the sequel, we denote by $DF(u)$ and
$D^2 F(u)$ the first and second Fr\'{e}chet derivatives of $F$ at the
point $u\in H_0^1(B)$.

\begin{lemma}
Assume $u$ is a nonzero critical point of $F$, normalized with $||u||=1$.
Then, for all $v\in H_0^1(B)$,
\begin{multline} \label{der-sec}
D^2F(u)[v,v]= \\
\gamma^2 \int_B e^{\gamma u^2} \left(
4u^2 v^2 + 4u^4 \left( \int_B \nabla u \cdot \nabla v \, dx\right)^2
-8u^3v \int_B \nabla u \cdot \nabla v \, dx \right)\, d\mu_\alpha \\
{}+\gamma \int_B e^{\gamma u^2} \left( 2v^2-2u^2\int_B |\nabla v|^2 \,
dx \right)\, d\mu_\alpha.
\end{multline}
\end{lemma}

\begin{proof}
Let $u$ be any nonzero critical point for $F$. Thus,
\begin{equation} \label{u-crit}
DF(u)[v]=0 \quad \hbox{for all $v\in H_0^1(B)$},
\end{equation}
where
\[
DF(u)[v]=\gamma \int_B e^{\gamma N(u)} DN(u)[v]\, d\mu_\alpha
\]
and
\[
DN(u)[v]=\frac{2uv\int_B |\nabla u|^2\, dx - 2u^2 \int_B \nabla u
\cdot \nabla v \, dx}{\left( \int_B |\nabla u|^2\, dx \right)^2}.
\]
For every $v,w\in H_0^1(B)$, the second derivative of $F$
at $u$ is
\begin{align}
\label{dedue}
D^2F(u)[v,w]= \gamma \int_B e^{\gamma N(u)} &DN(u)[v] \, DN(u)[w]\, d\mu_\alpha
\nonumber \\
&+\gamma\int_B e^{\gamma N(u)} D^2 N(u)[v,w]\, d\mu_\alpha.
\end{align}
We now compute the two integrals. We have
\begin{align*}
&D^2 N(u)[v,w]=\left( \int_B |\nabla u|^2 \, dx \right)^{-4} \left(
\left( \int_B |\nabla u|^2 \, dx \right)^{2} \left( 2vw \int_B |\nabla
u|^2\, dx \right. \right. \\ &\left. \left. {} + 4uv \int_B \nabla u
\cdot \nabla w \, dx -4 uw \int_B \nabla u \cdot \nabla v \, dx - 2u^2
\int_B \nabla v \cdot \nabla w\, dx \right) \right. \\ &\left. {} - 8
\left( uv \int_B |\nabla u|^2 \, dx - u^2 \int_B \nabla u \cdot \nabla
v \, dx \right) \cdot \int_B |\nabla u|^2 \, dx \cdot \int_B \nabla u
\cdot \nabla w \, dx \right) \\ &= \left( \int_B |\nabla u|^2 \, dx
\right)^{-2} \left( 2 vw\int_B |\nabla u|^2 \, dx +4uv\int_B \nabla u
\cdot \nabla w\, dx \right.\\ &\left. \quad {}-4uw\int_B \nabla u
\cdot \nabla v\, dx -2u^2 \int_B \nabla v \cdot \nabla w \, dx \right)
\\ &{} \quad - \frac{4 \int_B \nabla u \cdot \nabla w\,dx}{\left( \int_B
|\nabla u|^2 \, dx \right)^{2}} \;\; \frac{2uv \int_B |\nabla u|^2\, dx
- 2u^2 \int_B \nabla u \cdot \nabla v\, dx}{\left( \int_B |\nabla u|^2
\, dx \right)^{2}}.
\end{align*}
If we recall \eqref{u-crit}, we conclude that
\begin{multline*}
\gamma \int_B e^{\gamma N(u)} D^2 N(u)[v,w]\, d\mu_\alpha = \\
\frac{\gamma}{\left( \int_B |\nabla u|^2 \, dx \right)^2}
\int_B e^{\gamma N(u)}\left( 2vw \int_B |\nabla u|^2\, dx
+ 4uv \int_B \nabla u \cdot \nabla w \, dx \right. \\
\left. {} -4 uw \int_B \nabla u \cdot \nabla v \, dx
- 2u^2 \int_B \nabla v \cdot \nabla w\, dx \right)\, d\mu_\alpha.
\end{multline*}
Therefore, choosing $w=v$, we immediately see that
\begin{multline*}
\gamma \int_B e^{\gamma N(u)} D^2 N(u)[v,v]\, d\mu_\alpha = \\
\frac{\gamma}{\left( \int_B |\nabla u|^2 \, dx \right)^2}
\int_B e^{\gamma N(u)}\left( 2v^2 \int_B |\nabla u|^2\,dx
- 2u^2 \int_B |\nabla v|^2\,dx \right)\, d\mu_\alpha.
\end{multline*}
If in addition $u$ is normalized by $\int_B |\nabla u|^2 \, dx =1$, then
\begin{equation}
\label{sec}
\gamma \int_B e^{\gamma N(u)} D^2 N(u)[v,v]\, d\mu_\alpha =
\gamma \int_B e^{\gamma u^2} \left( 2v^2 -2u^2
\int_B |\nabla v|^2 \, dx \right)\, d\mu_\alpha.
\end{equation}
As far as the first integral in (\ref{dedue}) is concerned, by similar but simpler
arguments, we obtain, for a normalized critical point,
\begin{equation}
\label{fir}
\left(
DN(u)[v]\right)^2 = 4u^2v^2 + 4u^4 \left( \int_B \nabla u \cdot \nabla
v \,dx\right)^2 - 8u^3 v \int_B \nabla u \cdot \nabla v\,dx.
\end{equation}
Finally, equation \eqref{der-sec} is an immediate consequence of
(\ref{sec}) and(\ref{fir}).
\qed\end{proof}

We can now prove the main result of this paper.

\begin{proof}[of Theorem \ref{main}]
Let $u=u_\ge$ be any solution to problem \eqref{P}, and assume that it
is a radial function. Any $v\in H_0^1(B)$ can be decomposed as
$v=au+w$ with $a\in \R$ and $\int_B \nabla u \cdot \nabla w\, dx
=0$. It follows from \eqref{der-sec} that
\begin{multline*}
D^2F(u)[au+w,au+w]= \\ 4\gamma^2 \int_B e^{\gamma u^2} u^2w^2
|x|^\alpha\, dx + 2\gamma \int_B e^{\gamma u^2} \left( w^2-u^2 \int_B
|\nabla w|^2 \right) |x|^\alpha \, dx.
\end{multline*}
Choose now $w=u\psi f$, where $\psi$ is a radial function and
$f(\theta)=\sin \theta$ (in polar coordinates). Then, using the fact that
$$
\int_0^{2\pi} f\,d\theta = 0\qquad\hbox{and}\qquad
\int_0^{2\pi} f^2\,d\theta = \int_0^{2\pi} f_\theta^2\,d\theta,
$$
we see that
\begin{multline} \label{der-sec-fin}
D^2F(u)[au+u\psi f,au+u\psi f]= \\ 4\gamma^2 \int_B e^{\gamma u^2} u^4
\psi^2 |x|^\alpha\, dx + 2\gamma \int_B e^{\gamma u^2} u^2 \psi^2
|x|^\alpha \, dx \\ {} - 2\gamma \int_B e^{\gamma u^2} u^2 |x|^\alpha
\, dx \left[ \int_B |\nabla (u\psi)|^2\, dx + \int_B \frac{u^2
\psi^2}{|x|^2}\, dx \right].
\end{multline}
Since $u$ is a solution to \eqref{P}, $D^2F(u)$ must be negative
semidefinite as a bilinear form on $H_0^1(B)$. We choose a suitable
$\psi$ and deduce that this can hold (for $\alpha$ large) only if
\begin{equation}
\label{condit}
\gamma \leq \frac{\pi \varphi_1(0)^2}{\lambda_1 \int_B \varphi_1^4 \, dx}.
\end{equation}
We take $\psi (r)=r$, and refer to the last remark why we
choose this simple candidate.

By direct computation, if $\psi(r)=r$, then
\[
\int_B |\nabla (u\psi)|^2\, dx + \int_B \frac{u^2 \psi^2}{|x|^2}\, dx
= \int_B |\nabla u|^2 |x|^2\, dx,
\]
and therefore
\begin{multline} \label{eq:29}
D^2F(u)[au+u\psi f,au+u\psi f]= \\ 4\gamma^2 \int_B e^{\gamma u^2} u^4
|x|^{\alpha+2}\, dx + 2\gamma \int_B e^{\gamma u^2} u^2 |x|^{\alpha+2}
\, dx \\ {} - 2\gamma \int_B e^{\gamma u^2} u^2 |x|^\alpha \, dx
\int_B |\nabla u|^2 |x|^2\, dx.
\end{multline}
Recalling that $u$ satisfies the Euler--Lagrange equation
\eqref{eulero} with
$$
\lambda = \left(\int_B e^{\gamma u^2} u^2 |x|^\alpha \, dx
\right)^{-1},
$$
we can multiply both sides of equation \eqref{eulero} by $|x|^2 u$ and
integrate to obtain
\begin{equation} \label{eq:29bis}
2\gamma\lambda \int_B e^{\gamma u^2} u^2 |x|^{\alpha+2} \, dx = \int_B
|\nabla u|^2 |x|^2\, dx - 2\int_B u^2\, dx.
\end{equation}
We remark that we could have found the last identity by a direct use
of condition \eqref{u-crit}.  Inserting \eqref{eq:29bis} into
\eqref{eq:29} we find
\begin{multline*}
D^2F(u)[au+u\psi f,au+u\psi f]= \\ 4\gamma^2 \int_B e^{\gamma u^2} u^4
|x|^{\alpha+2}\, dx - 4\gamma \int_B e^{\gamma u^2} u^2 |x|^\alpha \,
dx \int_B u^2\, dx.
\end{multline*}
We write $u(|x|)=\sqrt{\ge} v_\ge (|x|^{1/\ge})$, and
we recall from Theorem~\ref{th2} that $v_\ge \to \varphi_1$
strongly in $H_0^1(B)$. Plugging the new variable $v_\ge$ into the
previous equation gives
\begin{multline} \label{eq:30}
D^2F(u)[au+u\psi f,au+u\psi f]= 4\gamma \ge^3 \Big\{ \gamma \int_B
e^{\gamma \ge v_\varepsilon^2} v_\varepsilon^4 |x|^{2\ge}\, dx \\
{} -\int_B e^{\gamma \ge
v_\varepsilon^2} v_\varepsilon^2 \, dx \cdot \ge \int_0^{2\pi} \int_0^1 v_\varepsilon^2
t^{2\ge -1} \, dt\, d\theta \Big\} .
\end{multline}
By a simple integration by parts, one checks immediately that
\[
\lim_{\ge\to 0} \ge \int_0^{2\pi} \int_0^1 v_\varepsilon^2
t^{2\ge -1} \, dt\, d\theta =\pi \varphi_1(0)^2.
\]
We may now conclude that
\begin{multline*}
\lim_{\ge \to 0} \frac{1}{4\gamma \ge^3} D^2F(u)[au+u\psi f,au+u\psi f] = \\
\gamma \int_B \varphi_1^4\, dx - \pi \varphi_1(0)^2\int_B \varphi_1^2
= \gamma \int_B \varphi_1^4\, dx - \frac{\pi}{\lambda_1}
\varphi_1(0)^2 \leq 0
\end{multline*}
only if condition (\ref{condit}) holds. This completes the proof.
\qed\end{proof}

\begin{remark}
We give a formal motivation why we have chosen $\psi (|x|)=|x|$ in the
proof of the theorem. It is clear that equation~\eqref{der-sec-fin} is
homogenous in $\psi$, so that we can assume without loss of generality
$\psi(1)=1$. By inspecting \eqref{der-sec-fin}, in order to make
$D^2F(u)[au+u\psi f,au+u\psi f]$ negative it seems natural to choose a
$\psi$ among radial functions vanishing at zero which keeps the
integral
\[
\int_B \Big( |\nabla \psi|^2\, +\frac{\psi^2}{|x|^2} \Big) u^2\, dx
\]
as small as possible. The heuristic reason why we have chosen
$\psi(r)=r$ is that this is precisely the unique solution to the
variational problem
\[
\inf \left\{ \int_B \Big( |\nabla \psi|^2\, +\frac{\psi^2}{|x|^2}
\Big)\, dx : \hbox{$\psi(1)=1$ and $\psi(0)=0$} \right\}.
\]
\end{remark}

\begin{remark}
Our results lead, in a natural way, to the following question: does
there exist a ``bifurcation point'' $\gamma_\star<4\pi$ such that
non-radial maximizers of problem \eqref{P} exist only when $\gamma >
\gamma_\star$?
\end{remark}

\begin{acknowledgement}
We would like to thank M.~Calanchi and E.~Terraneo (both University of
Milan) for some useful discussions about our problem.
\end{acknowledgement}



\end{document}